\input amstex

\define\Pee{{\Bbb P}}
\define\Zee{{\Bbb Z}}

\define\W{\widetilde}
\define\w{\tilde}

\define\Pic{\operatorname{Pic}}
\define\Num{\operatorname{Num}}

\define\Hilb{\operatorname{Hilb}}

\define\proof{\demo{Proof}}
\define\endproof{\qed\enddemo}

\define\theorem#1{\proclaim{Theorem #1}}
\define\lemma#1{\proclaim{Lemma #1}}
\define\proposition#1{\proclaim{Proposition #1}}
\define\corollary#1{\proclaim{Corollary #1}}
\define\claim#1{\proclaim{Claim #1}}

\define\section#1{\specialhead #1 \endspecialhead}
\define\ssection#1{\medskip\noindent{\bf #1}}

\loadbold
\documentstyle{amsppt}
\leftheadtext{Wei-ping Li and Zhenbo Qin}
\rightheadtext{Blowup formulae}

\topmatter
\title
On blowup formulae for the $S$-duality conjecture of Vafa and Witten II:
the universal functions
\endtitle
\author {Wei-ping Li$^1$ and Zhenbo Qin$^2$}
\endauthor
\address 
Department of Mathematics, HKUST, Clear Water Bay, Kowloon, Hong Kong
\endaddress
\email mawpli\@uxmail.ust.hk  \endemail
\address Department of Mathematics, Oklahoma State University, 
Stillwater, OK 74078, USA
\endaddress
\email  zq\@math.okstate.edu \endemail
\thanks
${}^1$Partially supported by the grant HKUST631/95P
\endthanks
\thanks 
${}^2$Partially supported by NSF grant DMS-9622564 and 
an Alfred P. Sloan Research Fellowship
\endthanks

\endtopmatter

\NoBlackBoxes
\TagsOnRight
\document
\section{1. Introduction}

This is a continuation of our work \cite{L-Q} on blowup formulae
for the S-duality conjecture of Vafa and Witten. 
In \cite{V-W}, Vafa and Witten formulated some mathematical predictions 
about the Euler characteristics of instanton moduli spaces
derived from the $S$-duality conjecture in physics. 
 From these mathematical predictions, a blowup formula was proposed
based upon the work of Yoshioka \cite{Yos}. 
Roughly speaking, the blowup formula says that there exists 
a universal relation between the Euler characteristics of 
instanton moduli spaces for a smooth four manifold
and the Euler characteristics of instanton moduli spaces for 
the blowup of the smooth four manifold. The universal relation is
independent of the four manifold and related to some modular forms.
In \cite{L-Q}, we verified this blowup formula for 
the gauge group $SU(2)$ and its dual group $SO(3)$ 
when the underlying four manifold is an algebraic surface.
In fact, we proved a stronger blowup formula in \cite{L-Q},
i.e. a blowup formula for the virtual Hodge numbers of instanton moduli spaces.
However, in \cite{L-Q}, we did not find a closed formula for 
the universal function which appears in this stronger blowup formula. 
Our goal of the present paper is to  determine a closed formula 
for this universal function. 

To state the blowup formulae proved in \cite{L-Q},
we recall some standard definitions and notations. 
Let $\phi: \W X \to X$ be the blowing-up of an algebraic surface $X$ 
at a point $x_0 \in X$, and $E$ be the exceptional divisor. 
For simplicity, we always assume that $X$ is simply connected. 
Fix a divisor $c_1$ on $X$, $\w c_1 = \phi^*c_1 - aE$ with $a = 0$ or $1$,
and an ample divisor $H$ on $X$ with odd $(H \cdot c_1)$. 
For an integer $n$, let $\frak M_{H}(c_1, n)$ be the moduli space
of Mumford-Takemoto $H$-stable rank-$2$ bundles with 
Chern classes $c_1$ and $n$, $\frak M^G_{H}(c_1, n)$ 
be the moduli space of Gieseker $H$-semistable rank-$2$ torsion-free sheaves 
with Chern classes $c_1$ and $n$, and $\frak M^U_{H}(c_1, n)$ be 
the Uhlenbeck compactification of $\frak M_{H}(c_1, n)$ from 
gauge theory \cite{Uhl, Don, LiJ}. It is well-known that 
both the Gieseker moduli spaces and the Uhlenbeck compactification spaces
are projective. For $r \gg 0$, the divisors $H_r = r \cdot \phi^*H - E$ on 
$\W X$ is ample; moreover, all the moduli spaces $\frak M_{H_r}(\w c_1, n)$ 
(resp. $\frak M^G_{H_r}(\w c_1, n)$, $\frak M^U_{H_r}(\w c_1, n)$) 
can be naturally identified. So we shall use 
$\frak M_{H_\infty}(\w c_1, n)$ (resp. $\frak M^G_{H_\infty}(\w c_1, n)$, 
$\frak M^U_{H_\infty}(\w c_1, n)$) to denote the moduli space
$\frak M_{H_r}(\w c_1, n)$ (resp. $\frak M^G_{H_r}(\w c_1, n)$, 
$\frak M^U_{H_r}(\w c_1, n)$) with $r \gg 0$. 

For a complex algebraic scheme $Y$ 
(not necessarily smooth, projective, or irreducible), 
let $e(Y; x, y)$ be the virtual Hodge polynomial of $Y$. 
When $Y$ is projective, $e(Y; 1,1)$ is the topological Euler characteristic 
of $Y_{\text{red}}$. Our Theorem A in \cite{L-Q} gives 
the following blowup formula for the Gieseker moduli spaces:
$$\sum_{n} e(\frak M_{H_\infty}^G(\w c_1, n); x, y) q^{n- {\w c_1^2 \over 4}} 
= (q^{1 \over 12} \cdot \W {\W Z}_a) \cdot 
\sum_{n} e(\frak M_{H}^G(c_1, n); x, y) q^{n- {c_1^2  \over 4}} \eqno (1.1)$$
where $\W {\W Z}_a = \W {\W Z}_a(x, y, q)$
is a universal function of $x, y, q, a$ with 
$$\W {\W Z}_a(1, 1, q) = {\sum_{n \in \Zee} q^{(n+{a \over 2})^2} \over 
[q^{1 \over 24} \prod_{n \ge 1} (1 - q^n)]^2}.$$
Assuming that $\frak M_H(c_1, n)$  
(respectively, $\frak M_{H_\infty}(\w c_1, n)$)
is dense in the Gieseker moduli space $\frak M_H^G(c_1, n)$
(respectively, $\frak M_{H_\infty}^G(\w c_1, n)$) for every $n$,
we also have a blowup formula for the Uhlenbeck compactification spaces
(Theorem B in \cite{L-Q}): 
$$\sum_{n} e(\frak M^U_{H_\infty}(\w c_1, n); x, y) 
q^{n- {\w c_1^2 \over 4}} = (q^{1 \over 12} \cdot \W Z_a) \cdot 
\sum_{n} e(\frak M^U_{H}(c_1, n); x, y) q^{n- {c_1^2 \over 4}}$$
where $\W Z_a = \W Z_a(x, y, q)$ is a universal function of $x, y, q, a$ with 
$$\W Z_a(1, 1, q) = {\sum_{n \in \Zee} q^{(n+{a \over 2})^2} \over 
q^{1 \over 12}(1 - q)}.$$

Our main results are the following closed formulae for 
$\W {\W Z}_a(x, y, q)$ and ${\W Z}_a(x, y, q)$.

\theorem{1.2} The universal function $\W {\W Z}_a(x, y, q)$ is equal to 
$${\sum_{n \in \Zee} (xy)^{(2n+a)^2 - (2n+a) \over 2} q^{(2n+a)^2 \over 4}
\over [q^{1 \over 24} \prod_{n \ge 1} (1 - (xy)^{2n}q^n)]^2}.$$
\endproclaim

\theorem{1.3} The universal function $\W Z_a(x, y, q)$ is equal to 
$$\align
&{1\over q^{1\over 12}(1-xyq)} \biggl [
\sum_{s\ge 0}(xy)^{(2s+a)^2+(2s+a) \over 2}q^{(2s+a)^2 \over 4}
\prod^{2s+a}_{j=1} {1-(xy)^{2j-2}q^j \over 1-(xy)^{2j}q^j}\\
&\qquad +\sum_{s\ge (1-a)}(xy)^{(2s+a)^2+(2s+a)-2 \over 2} q^{(2s+a)^2 \over 4}
\prod^{2s+a-1}_{j=1} {1-(xy)^{2j-2}q^j \over 1-(xy)^{2j}q^j} \biggr ] \\ 
\endalign$$
where we make the convention that 
$\prod \limits^{0}_{j=1} {1-(xy)^{2j-2}q^j \over 1-(xy)^{2j}q^j} = 1$.
\endproclaim

The paper is organized as follows. 
In section two, we verify Theorem 1.2 by taking $X = \Bbb F_1$
(the one-point blownup of $\Pee^2$). In section three, 
we prove Theorem 1.3 by using a not-closed formula of 
$\W Z_a(x, y, q)$ obtained in \cite{L-Q}.

Just before we post this paper, G\" ottsche informed us that 
in his forthcoming paper, he deals with the Hodge numbers of
the Gieseker moduli spaces on rational surfaces. 


\section{2. The universal function $\W {\W Z}_a(x, y, q)$}

In this section, we derive a closed formula for 
the universal function $\W {\W Z}_a(x, y, q)$.
Our strategy is to compute the virtual Hodge polynomials of
the Gieseker moduli spaces of semistable rank-$2$ sheaves over 
$\Bbb F_1$ and its blownup. These Gieseker moduli spaces are actually smooth 
and have been studied extensively (see \cite{E-G, F-Q} for example). 
Adopting a formula of G\" ottsche \cite{Got},
we calculate the (virtual) Hodge polynomials of these Gieseker moduli spaces.
Then using the definition of $\W {\W Z}_a(x, y, q)$, 
we can determine a closed formula for $\W {\W Z}_a(x, y, q)$.

First of all, we recall virtual Hodge polynomials for 
complex algebraic schemes
(not necessarily smooth, projective, or irreducible). 
Virtual Hodge polynomials were introduced 
by Danilov and Khovanskii \cite{D-K}. 
They can be viewed as a tool for computing the Hodge numbers of 
smooth projective varieties by reducing to 
computing those of simpler varieties. 
For a complex algebraic scheme $Y$, 
Deligne \cite{Del} proved that the cohomology $H_c^k(Y, \Bbb Q)$ 
with compact support carries a natural mixed Hodge structure 
which coincides with the classical one if $Y$ is projective and smooth. 
For each pair of integers $(s, t)$, define the virtual Hodge number
$$e^{s, t}(Y) = \sum_k (-1)^k h^{s, t}(H_c^k(Y, \Bbb Q)).$$ 
Then the virtual Hodge polynomials of $Y$ is defined by 
$$e(Y; x, y) = \sum_{s, t} e^{s, t}(Y) x^sy^t.  $$
Virtual Hodge polynomials satisfy the following properties 
(see \cite{D-K, Ful, Che}):

\roster
\item"{(2.1)}" When $Y$ is projective, $e(Y; 1, 1)$ is  
the Euler characteristic $\chi(Y_{\text{red}})$ of $Y_{\text{red}}$.
When $Y$ is projective and smooth, $e(Y; x, y)$ is the usual Hodge polynomial. 
\item"{(2.2)}" If $Z$ is a Zariski-closed subscheme of $Y$, then
$$e(Y; x, y) = e(Z; x, y) + e(Y - Z; x, y).$$
So if $Y = \coprod_{i=1}^n Y_i$ is a disjoint union of
finitely many locally closed subsets (i.e. each $Y_i$ is 
the intersection of an open subset and a closed subset), then
$$e(Y; x, y) = \sum_{i=1}^n e(Y_i; x, y).$$
\item"{(2.3)}" If $f: Y \to Z$ is a Zariski-locally trivial bundle 
with fiber $F$, then
$$e(Y; x, y) = e(Z; x, y) \cdot e(F; x, y).$$
\item"{(2.4)}" If $f: Y \to Z$ is a bijective morphism, 
then $e(Y; x, y) = e(Z; x, y)$. 
In particular, we have $e(Y; x, y) = e(Y_{\text{red}}; x, y)$.
\endroster

Next, we recall a result of G\" ottsche. Let $X$ be an algebraic surface with
effective anti-canonical divisor $-K_X$, and let $q(X)$ be its irregularity.
Fix a divisor $c_1$ on $X$ and an integer $n$. In \cite{Got}, 
G\" ottsche studied the change of the virtual Hodge polynomial 
$e(\frak M^G_H(c_1, n);x,y)$ as the ample divisor $H$ crosses walls of 
type $(c_1, n)$. In addition, a detailed study of the change of 
the Gieseker moduli space $\frak M^G_H(c_1, n)$ as $H$ crosses walls of 
type $(c_1, n)$ can be found in \cite{E-G, F-Q}.
The next lemma follows immediately from the Theorem 3.4 (1) in \cite{Got}.

\lemma{2.5} Assume that $X$ is an algebraic surface with
effective $-K_X$. Let $H$ and $L$ be ample divisors not lying
on any wall of type $(c_1, n)$. Then
$$e(\frak M^G_H(c_1, n);x,y) = e(\frak M^G_L(c_1, n);x,y)
+((1-x)(1-y))^{q(X)} \cdot  $$
$$\cdot \sum_\zeta (xy)^{\ell_\zeta - {\zeta^2 + \zeta K_X \over 2} -
\chi(\Cal O_X)} {1 -(xy)^{\zeta K_X} \over 1-(xy)} \cdot
\sum_{s+t=\ell_\zeta} e(\Hilb^s(X);x,y) e(\Hilb^t(X);x,y) $$
where $\ell_\zeta = (4n - c_1^2 + \zeta^2)/4$, and $\zeta$ rus over
all the classes in $\Num(X)$ which
define walls of type $(c_1, n)$ and satisfy $\zeta H < 0 < \zeta L$. \qed
\endproclaim

Now let $X$ be a rational ruled surface with effective $-K_X$.
Then $q(X) = 0$. Let $f$ be a general fiber of the ruling.
Fix a divisor $c_1$ and an ample divisor $H$ such that
both $(f \cdot c_1)$ and $(H \cdot c_1)$ are odd. Fix an integer $n$.
Since $(H \cdot c_1)$ is odd, $H$ does not lie on any wall of type $(c_1, n)$.
Since $(f \cdot c_1)$ is odd, it is well-known \cite{H-S, Qi2} that
there exists an open chamber $\Cal C_n$ of type $(c_1, n)$
such that $\frak M^G_{L_n}(c_1, n) = \emptyset$ for $L_n \in \Cal C_n$
and that the divisor class $f$ is contained in the closure of $\Cal C_n$.
Note that since the divisor $f$ is nef and contained in
the closure of $\Cal C_n$, the condition $\zeta H < 0 < \zeta L_n$ 
is equivalent to $\zeta H < 0 < \zeta f$. Let
$$\Lambda_H = \{ \zeta \in \Pic(X) | \quad \zeta H < 0 < \zeta f
\text{ and } \zeta \equiv c_1 \pmod 2 \}.  \eqno (2.6)$$
Then $\zeta$ defines a nonempty wall of type $(c_1, n)$ with
$\zeta H < 0 < \zeta L_n$ if and only if $\zeta \in \Lambda_H$
and $\zeta^2 \ge -(4n-c_1^2)$. Applying Lemma 2.5 to $H$ and $L_n$,
we obtain
$$e(\frak M^G_H(c_1, n);x,y) = 
\sum_{\zeta \in \Lambda_H \text{ and } \zeta^2 \ge -(4n-c_1^2)}
(xy)^{\ell_\zeta - {\zeta^2 + \zeta K_X \over 2} -
\chi(\Cal O_X)} {1 -(xy)^{\zeta K_X} \over 1-(xy)} \cdot$$
$$\cdot \sum_{s+t=\ell_\zeta} e(\Hilb^s(X);x,y) e(\Hilb^t(X);x,y).
\eqno (2.7)$$

\lemma{2.8} Let $X$ be a rational ruled surface with effective $-K_X$. 
Let $c_1$ be a divisor on $X$ such that both $(f \cdot c_1)$ and 
$(H \cdot c_1)$ are odd. Then 
$$\sum_n e(\frak M^G_H(c_1, n);x,y) q^{n - {c_1^2 \over 4}} =
{[\sum_n e(\Hilb^n(X);x,y) (xyq)^n]^2 \over
(xy)^{\chi(\Cal O_X)}[1-(xy)]} \cdot$$
$$\cdot \sum_{\zeta \in \Lambda_H} (xy)^{-{\zeta^2 + \zeta K_X \over 2}}
[1 -(xy)^{\zeta K_X}] q^{-{\zeta^2 \over 4}}.  \eqno (2.9)$$
\endproclaim
\noindent
{\it Proof.} By definition, $\ell_\zeta = (4n - c_1^2 + \zeta^2)/4 \ge 0$.
So $n = \ell_\zeta + (c_1^2 - \zeta^2)/4$. By (2.7), 
$$\align
&\qquad \sum_n e(\frak M^G_H(c_1, n);x,y) q^{n - {c_1^2 \over 4}}  \\
&= \sum_n \quad
   \sum_{\zeta \in \Lambda_H \text{ and } \zeta^2 \ge -(4n-c_1^2)}
   (xy)^{\ell_\zeta - {\zeta^2 + \zeta K_X \over 2} -
   \chi(\Cal O_X)} {1 -(xy)^{\zeta K_X} \over 1-(xy)} \cdot \\
&\qquad \cdot \sum_{s+t=\ell_\zeta} e(\Hilb^s(X);x,y) e(\Hilb^t(X);x,y)
q^{n - {c_1^2 \over 4}} \\
&= \sum_{\zeta \in \Lambda_H}  \sum_{\ell \ge 0}
   (xy)^{\ell - {\zeta^2 + \zeta K_X \over 2} -
   \chi(\Cal O_X)} {1 -(xy)^{\zeta K_X} \over 1-(xy)} \cdot \\
&\qquad \cdot \sum_{s+t=\ell} e(\Hilb^s(X);x,y) e(\Hilb^t(X);x,y)
q^{\ell - {\zeta^2 \over 4}} \\
&= \sum_{\zeta \in \Lambda_H} (xy)^{-{\zeta^2 + \zeta K_X \over 2} -
   \chi(\Cal O_X)} {1 -(xy)^{\zeta K_X} \over 1-(xy)}
   q^{-{\zeta^2 \over 4}} \cdot \\
&\qquad \cdot \sum_{\ell \ge 0} \sum_{s+t=\ell} e(\Hilb^s(X);x,y)
e(\Hilb^t(X);x,y) (xyq)^{\ell}. \\
\endalign$$
Here going from the first equality to the second equality, 
we have changed $n$ to $\ell + (c_1^2 - \zeta^2)/4$ with $\ell \ge 0$. 
Notice that 
$$\sum_{\ell \ge 0} \sum_{s+t=\ell} e(\Hilb^s(X);x,y)
e(\Hilb^t(X);x,y) (xyq)^{\ell}$$ 
is equal to $[ \sum_n e(\Hilb^n(X);x,y) (xyq)^n ]^2$. Therefore, we obtain
$$\align
&\qquad \sum_n e(\frak M^G_H(c_1, n);x,y) q^{n - {c_1^2 \over 4}}  \\
&= \sum_{\zeta \in \Lambda_H} (xy)^{-{\zeta^2 + \zeta K_X \over 2} -
   \chi(\Cal O_X)} {1 -(xy)^{\zeta K_X} \over 1-(xy)}
   q^{-{\zeta^2 \over 4}} \cdot 
   \left [\sum_n e(\Hilb^n(X);x,y) (xyq)^n \right ]^2  \\
&= {[\sum_n e(\Hilb^n(X);x,y) (xyq)^n]^2 \over
(xy)^{\chi(\Cal O_X)}[1-(xy)]} \cdot 
\sum_{\zeta \in \Lambda_H} (xy)^{-{\zeta^2 + \zeta K_X \over 2}}
[1 -(xy)^{\zeta K_X}] q^{-{\zeta^2 \over 4}}. \qed \\
\endalign$$

Next we study the virtual Hodge polynomials of the Gieseker moduli spaces
over blownup surfaces. As before, let $X$ be a rational ruled surface 
with effective $-K_X$. Let $f$ be a general fiber of the ruling.
Fix a divisor $c_1$ and an ample divisor $H$ on $X$ such that
both $(f \cdot c_1)$ and $(H \cdot c_1)$ are odd. 
Let $\phi: \W X \to X$ be the blowing-up of $X$ at a point $x_0 \in X$, 
and $E$ be the exceptional divisor. We assume that $-K_{\W X}$ is effective.
Let $\w c_1 = \phi^*c_1 - aE$ with $a = 0$ or $1$. 
It is well-known \cite{F-M, Bru, Qi1} that for $r \gg 0$, 
all the divisors $H_r = r \cdot \phi^*H - E$ on $\W X$ are ample and  
lie in the same open chamber of type $(\w c_1, n)$. Thus all the moduli spaces 
$\frak M^G_{H_r}(\w c_1, n)$ (resp. $\frak M_{H_r}(\w c_1, n)$) 
with $r \gg 0$ are identical, and shall be denoted by 
$\frak M^G_{H_\infty}(\w c_1, n)$ (resp. $\frak M_{H_\infty}(\w c_1, n)$).
Since $(H_r \cdot \w c_1) = r(H \cdot c_1) - a$ and $(H \cdot c_1)$ is odd,
we can always choose $r \gg 0$ such that $(H_r \cdot \w c_1)$ is also odd.

\lemma{2.10} Let $\phi: \W X \to X$ be the blowing-up of 
a rational ruled surface $X$ at one point such that 
$-K_X$ and $-K_{\W X}$ are effective. Let $c_1$ be a divisor on $X$ 
such that both $(f \cdot c_1)$ and $(H \cdot c_1)$ are odd,
and $\w c_1 = \phi^*c_1 - aE$ with $a = 0$ or $1$. Then 
$$\sum_n e(\frak M^G_{H_\infty}(\w c_1, n);x,y) q^{n - {\w c_1^2 \over 4}} =
{[\sum_n e(\Hilb^n(\W X);x,y) (xyq)^n]^2 \over
(xy)^{\chi(\Cal O_{\W X})}[1-(xy)]} \cdot$$
$$\cdot \sum_{t \in \Zee} (xy)^{(2t+a)^2 - (2t+a) \over 2} q^{(2t+a)^2 \over 4}
\cdot \sum_{\zeta \in \Lambda_H} (xy)^{-{\zeta^2 + \zeta K_X \over 2}}
[1 -(xy)^{\zeta K_X}] q^{-{\zeta^2 \over 4}}.  \eqno (2.11)$$
\endproclaim
\proof
Note that the ruling of $X$ induces a ruling of $\W X$ and 
that $\phi^*f$ is the divisor class of a general fiber for 
the ruling of $\W X$. Fix an integer $n$, 
and choose $r \gg 0$ such that $(H_r \cdot \w c_1)$ is odd. 
Applying (2.7) to $\W X$ and $H_r$, we obtain
$$e(\frak M^G_{H_r}(\w c_1, n);x,y) = 
\sum_{\w \zeta \in \Lambda_{H_r} \text{ and } {\w \zeta}^2 \ge -(4n-\w c_1^2)}
(xy)^{\ell_{\w \zeta} - {{\w \zeta}^2 + {\w \zeta} K_{\W X} \over 2} -
\chi(\Cal O_{\W X})} {1 -(xy)^{\w \zeta K_{\W X}} \over 1-(xy)} \cdot$$
$$\cdot \sum_{s+t=\ell_{\w \zeta}} e(\Hilb^s(\W X);x,y) e(\Hilb^t(\W X);x,y)
\eqno (2.12)$$
where by (2.6), $\Lambda_{H_r} = \{ \w \zeta \in \Pic(\W X) | 
\quad \w \zeta H_r < 0 < \w \zeta \cdot \phi^*f
\text{ and } \w \zeta \equiv \w c_1 \pmod 2 \}$.
Since $(H \cdot c_1)$ and $(H_r \cdot \w c_1)$ are odd,
$\phi^*H$ and $H_r$ are not separated by any wall of type $(\w c_1, n)$.
Thus if $\w \zeta$ defines a nonempty wall of type $(\w c_1, n)$,
then $\w \zeta H_r < 0 < \w \zeta \cdot \phi^*f$ if and only if 
$\w \zeta \cdot \phi^*H < 0 < \w \zeta \cdot \phi^*f$.
In view of this observation, we put
$$\Lambda_{H_\infty} = \{ \w \zeta \in \Pic(\W X) | 
\quad \w \zeta \cdot \phi^*H < 0 < \w \zeta \cdot \phi^*f
\text{ and } \w \zeta \equiv \w c_1 \pmod 2 \}.$$
Then by (2.12) and the convention for $\frak M^G_{H_\infty}(\w c_1, n)$, 
we have 
$$\align
&\qquad e(\frak M^G_{H_\infty}(\w c_1, n);x,y) 
= e(\frak M^G_{H_r}(\w c_1, n);x,y) = \\ 
&= \sum_{\w \zeta \in \Lambda_{H_\infty} \text{ and } 
  {\w \zeta}^2 \ge -(4n-\w c_1^2)} (xy)^{\ell_{\w \zeta} - {{\w \zeta}^2 + 
  {\w \zeta} K_{\W X} \over 2} - \chi(\Cal O_{\W X})} 
  {1 -(xy)^{\w \zeta K_{\W X}} \over 1-(xy)} \cdot\\
&\qquad \cdot \sum_{s+t=\ell_{\w \zeta}} e(\Hilb^s(\W X);x,y) 
  e(\Hilb^t(\W X);x,y). \\
\endalign$$
As in the proof of Lemma 2.8, we conclude that 
$$\sum_n e(\frak M^G_{H_\infty}(\w c_1, n);x,y) q^{n - {\w c_1^2 \over 4}} =
{[\sum_n e(\Hilb^n(\W X);x,y) (xyq)^n]^2 \over
(xy)^{\chi(\Cal O_{\W X})}[1-(xy)]} \cdot$$
$$\cdot \sum_{\w \zeta \in \Lambda_{H_\infty}} (xy)^{-{{\w \zeta}^2 + 
\w \zeta K_{\W X} \over 2}} [1 -(xy)^{\w \zeta K_{\W X}}] 
q^{-{\w \zeta^2 \over 4}}.  \eqno (2.13)$$
Put $\w \zeta = \phi^*\zeta + s E$. 
Then $\w \zeta \cdot \phi^*H < 0 < \w \zeta \cdot \phi^*f$ if and only if
$\zeta H < 0 < \zeta f$. Moreover, $\w \zeta \equiv \w c_1 \pmod 2$
if and only if $\zeta \equiv c_1 \pmod 2$ and $s \equiv a \pmod 2$.
So $\w \zeta = \phi^*\zeta + s E \in \Lambda_{H_\infty}$
if and only if $\zeta \in \Lambda_H$ and $s = (2t - a)$ for some $t \in \Zee$. Thus,
$$\align
&\qquad \sum_{\w \zeta \in \Lambda_{H_\infty}} (xy)^{-{{\w \zeta}^2 + 
\w \zeta K_{\W X} \over 2}} [1 -(xy)^{\w \zeta K_{\W X}}] 
q^{-{\w \zeta^2 \over 4}}\\
&= \sum_{\zeta \in \Lambda_H} \sum_{t \in \Zee} (xy)^{-{\zeta^2 - (2t-a)^2+ 
\zeta K_X - (2t-a) \over 2}} [1 -(xy)^{\zeta K_X - (2t-a)}] 
q^{-{\zeta^2 - (2t-a)^2 \over 4}} \\
&= \sum_{\zeta \in \Lambda_H} (xy)^{-{\zeta^2 + \zeta K_X \over 2}} 
q^{-{\zeta^2 \over 4}} \cdot \sum_{t \in \Zee} 
\left [ (xy)^{(2t-a)^2 + (2t-a) \over 2} - (xy)^{\zeta K_X + 
(2t-a)^2 - (2t-a) \over 2}  \right ] q^{(2t-a)^2 \over 4} \\
&= \sum_{\zeta \in \Lambda_H} (xy)^{-{\zeta^2 + \zeta K_X \over 2}} 
q^{-{\zeta^2 \over 4}} \cdot \sum_{t \in \Zee} 
\left [ (xy)^{(2t+a)^2 - (2t+a) \over 2} - (xy)^{\zeta K_X + 
(2t+a)^2 - (2t+a) \over 2}  \right ] q^{(2t+a)^2 \over 4} \\
&= \sum_{\zeta \in \Lambda_H} (xy)^{-{\zeta^2 + \zeta K_X \over 2}} 
[1 - (xy)^{\zeta K_X}] q^{-{\zeta^2 \over 4}} \cdot \sum_{t \in \Zee} 
(xy)^{(2t+a)^2 - (2t+a) \over 2} q^{(2t+a)^2 \over 4}. \tag 2.14 \\
\endalign$$
Here going from the second equality to the third equality, 
we have changed $t$ to $-t$ in the first term in the brackets
and $t$ to $t+a$ in the second term in the brackets.
Now the formula (2.11) follows from (2.13) and (2.14).
\endproof

\theorem{2.15} The universal function $\W {\W Z}_a(x, y, q)$ is equal to 
$${\sum_{n \in \Zee} (xy)^{(2n+a)^2 - (2n+a) \over 2} q^{(2n+a)^2 \over 4}
\over [q^{1 \over 24} \prod_{n \ge 1} (1 - (xy)^{2n}q^n)]^2}.$$
\endproclaim
\noindent
{\it Proof.} First of all, we notice from \cite{G-S} that 
for any algebraic surface $X$,
$$\sum_n e(\Hilb^n(X);x,y) q^n = \prod_{n \ge 1} 
\prod_{s, t = 0}^2 (1 - x^{s+n-1}y^{t+n-1}q^n)^{(-1)^{s+t+1}h^{s, t}(X)}
\eqno (2.16)$$
where $h^{s, t}(X)$ stands for the Hodge numbers of $X$.
Next, let $X = \Bbb F_1$ be the blownup of $\Pee^2$ at one point,
and let $\sigma$ be the exceptional divisor in $X$. 
Then $X$ is a ruled surface with effective $-K_X$.
Let $f$ be a fiber of the ruling.
Let $\phi: \W X \to X$ be the blowing-up of $X$ at one point. 
Then $-K_{\W X}$ is also effective. 
Let $H = \sigma + 2f$ and $c_1 = \sigma$. 
Then $(H \cdot c_1) = 1 = (f \cdot c_1)$. 
So $(H \cdot c_1)$ and $(f \cdot c_1)$ are odd.
Therefore the conditions in Lemma 2.8 and Lemma 2.10 are satisfied.
Note that $\chi(\Cal O_{\W X}) = \chi(\Cal O_X)$, 
$h^{s, t}(\W X) = h^{s, t}(X)$ when $(s, t) \ne (1, 1)$, 
and $h^{1,1}(\W X) = 1+h^{1,1}(X)$. By (2.16), 
$${\sum_n e(\Hilb^n(\W X);x,y) (xyq)^n \over \sum_n e(\Hilb^n(X);x,y) (xyq)^n}
= {1 \over \prod_{n \ge 1} (1 - (xy)^{2n}q^n)}.  \eqno (2.17)$$
Combining (2.9), (2.11), (2.17) with (1.1), we see that
$$\align
&\qquad \W {\W Z}_a(x, y, q)    \\
&= {1 \over q^{1 \over 12}} \cdot 
  {\sum_n e(\frak M^G_{H_\infty}(\w c_1, n);x,y) q^{n - {\w c_1^2 \over 4}} 
   \over \sum_n e(\frak M^G_H(c_1, n);x,y) q^{n - {c_1^2 \over 4}}}    \\
&= {1 \over q^{1 \over 12}} \cdot
  {[\sum_n e(\Hilb^n(\W X);x,y) (xyq)^n]^2 \over
  [\sum_n e(\Hilb^n(X);x,y) (xyq)^n]^2} \cdot 
  \sum_{t \in \Zee} (xy)^{(2t+a)^2 - (2t+a) \over 2} q^{(2t+a)^2 \over 4} \\
&= {\sum_{n \in \Zee} (xy)^{(2n+a)^2 - (2n+a) \over 2} q^{(2n+a)^2 \over 4}
\over [q^{1 \over 24} \prod_{n \ge 1} (1 - (xy)^{2n}q^n)]^2}.  \qed  \\
\endalign$$

\section{3. The universal function ${\W Z}_a(x, y, q)$}

In this section, we prove a closed formula for 
the universal function ${\W Z}_a(x, y, q)$. 
Our first goal is to compute the virtual Hodge polynomial of 
the space $U(m_1,m_2)$ which parameterizes all surjective maps 
$\Cal O_{\Pee^1}(-m_1) \oplus \Cal O_{\Pee^1}(-m_2) \to \Cal O_{\Pee^1} \to 0$. 
Then using the results in \cite{L-Q}, 
we obtain a closed formula for ${\W Z}_a(x, y, q)$. 
We end this section with a remark about this closed formula.

First of all, for two integers $m_1, m_2 \ge 0$, 
let $U(m_1, m_2)$ be the subset of 
$$\Pee(H^0(\Pee^1, \Cal O_{\Pee^1}(m_1) \oplus \Cal O_{\Pee^1}(m_2))) 
\cong \Pee^{m_1+m_2+1}$$
parameterizing all pairs $(f_1, f_2)$ of homogeneous polynomials 
such that $\deg (f_1) = m_1, \deg (f_2) = m_2$,
and $f_1$ and $f_2$ are coprime. Then $U(m_1, m_2)$ parameterizes 
all surjective maps $\Cal O_{\Pee^1}(-m_1) \oplus \Cal O_{\Pee^1}(-m_2)
\to \Cal O_{\Pee^1} \to 0$. The following result gives 
the virtual Hodge polynomial of $U(m_1, m_2)$.

\lemma{3.1} Let $m_1$ and $m_2$ be two integers with 
$0 \le m_1 \le m_2$. Then,
$$e(U(m_1, m_2); x, y) =  
\cases (xy)+1,                      &\text{if $m_1 = m_2 = 0$}        \\
       (xy)^{m_2+1},                &\text{if $m_1 = 0$ and $m_2 > 0$}\\
       (xy)^{m_1+m_2-1}[(xy)^2-1],  &\text{if $m_1 > 0$}.\\
\endcases  \eqno (3.2)$$
\endproclaim
\noindent
{\it Proof.}
We computed $e(U(m_1, m_2); 1, 1)$ in the Lemma 4.13 of \cite{L-Q}.
We shall adopt the same approach. First of all, 
we prove that (3.2) is true for $m_1 = 0$. Indeed, the subset $U(0, 0)$ of 
$\Pee(H^0(\Pee^1, \Cal O_{\Pee^1}(0) \oplus \Cal O_{\Pee^1}(0))) 
\cong \Pee^{1}$ coincides with $\Pee^1$. 
Since $e(\Pee^d; x, y) = 1 + (xy) + \ldots + (xy)^d$, we have  
$$e(U(0, 0); x, y) = e(\Pee^1; x, y) = (xy)+1.$$ 
So (3.2) holds for $m_1 = m_2 = 0$. When $m_2 > 0$, the subset $U(0, m_2)$ of 
$$\Pee(H^0(\Pee^1,\Cal O_{\Pee^1}(0) \oplus \Cal O_{\Pee^1}(m_2)))$$
is $\Pee(H^0(\Pee^1,\Cal O_{\Pee^1}(0) \oplus \Cal O_{\Pee^1}(m_2)))-
\Pee(\{ 0 \} \oplus H^0(\Pee^1,\Cal O_{\Pee^1}(m_2))) \cong \Pee^{m_2 + 1} - \Pee^{m_2}$.
Thus,
$$e(U(0, m_2); x, y) = e(\Pee^{m_2 + 1}; x, y) - e(\Pee^{m_2}; x, y)
= (xy)^{m_2+1}.$$ 
Hence (3.2) also holds for $m_1 = 0$ and $m_2 > 0$.

Next let $m_1 > 0$. The possible degree of the greatest common divisor
of a pair 
$$(f_1, f_2) \in \Pee(H^0(\Pee^1, \Cal O_{\Pee^1}(m_1) \oplus 
\Cal O_{\Pee^1}(m_2))) - 
\Pee(\{ 0 \} \oplus H^0(\Pee^1,\Cal O_{\Pee^1}(m_2)))$$
can be $0, \ldots, m_1$. For $d = 0, \ldots, m_1$, let $Y_d$ be the subset of 
$$\Pee(H^0(\Pee^1,\Cal O_{\Pee^1}(m_1) \oplus \Cal O_{\Pee^1}(m_2))) 
- \Pee(\{ 0 \} \oplus H^0(\Pee^1,\Cal O_{\Pee^1}(m_2)))$$ 
parameterizing all pairs $(f_1, f_2)$ such that 
$\hbox{gcd}(f_1, f_2)$ has degree $d$. Then we obtain
$$\align
\Pee^{m_1 +m_2+1} - \Pee^{m_2} 
& \cong \Pee(H^0(\Pee^1,\Cal O_{\Pee^1}(m_1) \oplus \Cal O_{\Pee^1}(m_2))) 
- \Pee(\{ 0 \} \oplus H^0(\Pee^1,\Cal O_{\Pee^1}(m_2))) \\
&= \coprod_{d = 0, \ldots, m_1} Y_d. \tag 3.3 
\endalign$$
Let $1 \le d \le m_1$, and $(f_1, f_2) \in Y_d$ with  
$\hbox{gcd}(f_1, f_2) = f$. Then we can write $f_1 = f g_1$ and $f_2 = f g_2$ 
with $f \in \Pee(H^0(\Pee^1,\Cal O_{\Pee^1}(d))) \cong \Pee^d$ and 
$$(g_1, g_2) \in 
\cases U(m_1-d, m_2-d), &\text{if $1 \le d < m_1$}\\
       U(0, m_2-m_1), &\text{if $d = m_1 < m_2$}\\ 
       U(0, 0) - \{ \text{a point} \}, &\text{if $d = m_1 = m_2$}.\\
\endcases$$
Thus $Y_d$ is the product of  the space $\Pee^d$ with the space 
 $U(m_1-d, m_2-d)$ when $1 \le d < m_1$ or $d = m_1 < m_2$,
or with the space $U(0, 0) - \{ \text{a point} \} \cong \Pee^1 - \{ \text{a point} \}$ 
when $d = m_1 = m_2$. So for $1 \le d \le m_1$, we have
$$e(Y_d; x, y) = e(\Pee^d; x, y) \cdot 
\cases  
e(U(m_1-d, m_2-d); x, y), &\text{if $1 \le d < m_1$}\\
e(U(0, m_2-m_1); x, y), &\text{if $d = m_1 < m_2$}\\
(xy), &\text{if $d = m_1 = m_2$.}\\
\endcases$$
Since $e(U(0, m_2-m_1); x, y) = (xy)^{m_2-m_1+1}$ when $m_1 < m_2$, we obtain
$$e(Y_d; x, y) = \sum_{i=0}^d (xy)^i \cdot 
\cases  
e(U(m_1-d, m_2-d); x, y), &\text{if $1 \le d < m_1$}\\
(xy)^{m_2-m_1+1}, &\text{if $d = m_1$.}\\
\endcases \tag 3.4$$
Note that $Y_0 = U(m_1, m_2)$. From (3.3) and (3.4), we conclude that
$$\align
\sum_{i=m_2+1}^{m_1+m_2+1} (xy)^i 
&= e(U(m_1, m_2); x, y) + \sum_{i=0}^{m_1} (xy)^i \cdot (xy)^{m_2-m_1+1} \\
&\quad+ \sum_{1 \le d < m_1} \sum_{i=0}^d (xy)^i 
 \cdot e(U(m_1-d, m_2-d); x, y). \tag 3.5 
\endalign$$

Now we see from (3.5) that $e(U(1, m_2); x, y) = (xy)^{m_2}[(xy)^2-1]$.
So (3.2) holds for $m_1 = 1$. For $m_1 > 1$, 
we use (3.5) and induction on $m_1$:
$$\align
e(U(m_1, m_2); x, y)
&= \sum_{i=m_2+1}^{m_1+m_2+1} (xy)^i - 
   \sum_{i=0}^{m_1} (xy)^i \cdot (xy)^{m_2-m_1+1} \\
&\quad- \sum_{1 \le d < m_1} \sum_{i=0}^d (xy)^i 
 \cdot (xy)^{(m_1-d)+(m_2-d)-1} [(xy)^2-1].\\
&= (xy)^{m_1+m_2-1} [(xy)^2-1].  \qed
\endalign$$

In section four of \cite{L-Q}, we proved the following formula:
$${\W Z}_a(x, y, q) = {q^{a \over 4} \cdot \sum_{n \ge 0} 
B_{a, n}(x,y) q^n \over q^{1 \over 12}(1 - xyq)}  \eqno (3.6)$$
where $B_{0, 0}(x,y)=1$, and $B_{a, n}(x,y)$ with $n \ge (1-a)$ is given by 
$$\align
B_{0, n}(x, y) =
&\sum_{\Sb 0 \le d_1,  0 \le d_{2j} \le d_{2j-1} -1,  
 0 \le d_{2j+1} \le d_{2j} (1 \le j \le s-1), 0 \le d_{2s} \le d_{2s-1} -1\\
 \sum_{i = 1}^{2s} d_i = n \endSb} \\
&\biggl (\prod_{i=1}^{s-1} e(U(d_{2i-1}-d_{2i} -1, d_{2i-1}+d_{2i}); x, y)  \tag 3.7\\
&e(U(d_{2i}-d_{2i+1}, d_{2i}+d_{2i+1}); x, y) \biggr )    \\
&e(U(d_{2s-1}-d_{2s} -1, d_{2s-1}+d_{2s}); x, y) e(U(d_{2s}, d_{2s}); x, y)  \\
B_{1, n}(x, y)= 
&\sum_{\Sb 0 \le d_1,  0 \le d_{2i} \le d_{2i-1}, 0 \le d_{2i+1} \le d_{2i} -1
       (1 \le i \le s)\\ \sum_{i = 1}^{2s+1} d_i = n \endSb} \\
&\biggl ( \prod_{i=1}^{s} e(U(d_{2i-1}-d_{2i}, d_{2i-1}+d_{2i}); x, y)  \tag 3.8\\
&e(U(d_{2i}-d_{2i+1}-1, d_{2i}+d_{2i+1}); x, y) \biggr ) 
e(U(d_{2s+1}, d_{2s+1}); x, y).    
\endalign$$

Now we can prove a closed formula for $\W Z_a(x, y, q)$.

\theorem{3.9} The universal function $\W Z_a(x, y, q)$ is equal to 
$$\align
&{1\over q^{1\over 12}(1-xyq)} \biggl [
\sum_{s\ge 0}(xy)^{(2s+a)^2+(2s+a) \over 2}q^{(2s+a)^2 \over 4}
\prod^{2s+a}_{j=1} {1-(xy)^{2j-2}q^j \over 1-(xy)^{2j}q^j}\\
&\qquad +\sum_{s\ge (1-a)}(xy)^{(2s+a)^2+(2s+a)-2 \over 2} q^{(2s+a)^2 \over 4}
\prod^{2s+a-1}_{j=1} {1-(xy)^{2j-2}q^j \over 1-(xy)^{2j}q^j} 
\biggr ] \tag 3.10 \\ 
\endalign$$
where we make the convention that 
$\prod \limits^{0}_{j=1} {1-(xy)^{2j-2}q^j \over 1-(xy)^{2j}q^j} = 1$.
\endproclaim
\noindent
{\it Proof.} Since the proof for the case $a=1$ is similar, 
we shall only prove the case $a=0$. By (3.6), it suffices to show that 
$$\align
\sum_{n \ge 0} B_{0, n}(x, y) q^n 
&= \sum_{s\ge 0} (xy)^{2s^2+s}q^{s^2} 
   \prod^{2s}_{j=1} {1-(xy)^{2j-2}q^j \over 1-(xy)^{2j}q^j}\\
&\quad +\sum_{s\ge 1}(xy)^{2s^2+s-1} q^{s^2}
\prod^{2s-1}_{j=1} {1-(xy)^{2j-2}q^j \over 1-(xy)^{2j}q^j}. \tag 3.11 \\ 
\endalign$$

First of all, let $\{d_1, d_2, \ldots, d_{2s} \}$ be an indexing sequence 
in the summation (3.7). So $0 \le d_1,  0 \le d_{2j} \le d_{2j-1} -1,  
 0 \le d_{2j+1} \le d_{2j} (1 \le j \le s-1), 0 \le d_{2s} \le d_{2s-1} -1$,
and $\sum \limits_{i = 1}^{2s} d_i = n$. 
We make the following chang of indices:
$$\cases  d_1' &=d_1-d_2-1,             \\
          d_2' &=d_2-d_3,               \\             
               &\vdots                  \\
     d_{2s-3}' &= d_{2s-3}-d_{2s-2}-1,  \\
     d_{2s-2}' &= d_{2s-2}-d_{2s-1},    \\
     d_{2s-1}' &=d_{2s-1}-d_{2s}-1,     \\
       d_{2s}' &=d_{2s}.
   \endcases         \eqno (3.12)$$
Thus, $d_i' \ge 0$ for all the $i$ with $1 \le i \le 2s$. Moreover, we have 
$$\cases    d_1 &=d_1'+\ldots +d_{2s}'+s,                 \\
            d_2 &=d_2'+\ldots +d_{2s}'+(s-1),             \\
            d_3 &=d_3'+\ldots +d_{2s}'+(s-1),             \\
                &\vdots                                   \\
       d_{2s-2} &=d_{2s-2}'+d_{2s-1}'+d_{2s}'+1,          \\
       d_{2s-1} &=d_{2s-1}'+d_{2s}'+1,                    \\
       d_{2s}   &=d_{2s}'.
\endcases \eqno(3.13)$$
So the condition $\sum\limits_{i = 1}^{2s} d_i=n$ becomes
$\sum \limits_{i=1}^{2s}id_i'+s^2 =n$ .

Next, let $t = (xy)$, and let $f: \{ 0, 1, 2, \ldots \} \to \{0, 1\}$ 
be defined by $f(d')= 0$ if $d'=0$, and $f(d')= 1$ if $d' > 0$.
Then for $0 \le d' \le d''$, (3.2) can be rewritten as:
$$e(U(d', d''); x, y) =  
\displaystyle{\cases t+1,                      &\text{if $d''= 0$}.\\
       t^{d'+d''+1}(1-\displaystyle{1\over t^2})^{f(d')}, &\text{if $d'' > 0$}.\\
\endcases}  \eqno (3.14)$$
Thus by (3.12), (3.13) and (3.14), the typical term in (3.7) is
$$\align
&\biggl (\prod_{i=1}^{s-1} e(U(d_{2i-1}-d_{2i} -1, d_{2i-1}+d_{2i}); x, y) 
 e(U(d_{2i}-d_{2i+1}, d_{2i}+d_{2i+1}); x, y) \biggr )  \\
&\quad e(U(d_{2s-1}-d_{2s} -1, d_{2s-1}+d_{2s}); x, y) e(U(d_{2s}, d_{2s}); x, y)\\
=& \cases t^{2d_1'+\ldots + 2d_{2s}' +2s} (1-{1\over t^2})^{f(d_1')}\ldots
          t^{2d_{2s}'+1}(1-{1\over t^2})^{f(d_{2s}')}, &\text{if $d_{2s}' > 0$,} \\
          t^{2d_1'+\ldots + 2d_{2s}' +2s }(1-{1\over t^2})^{f(d_1')}\ldots
          t^{2d_{2s-1}'+2}(1-{1\over t^2})^{f(d_{2s-1}')}(1+t), &\text{if $d_{2s}' = 0$.}
   \endcases \\
\endalign$$
It follows from (3.7) that $\sum_{n\ge 0}B_{0, n}(x, y)q^n$ is equal to 
$$\align
&1+\sum_{n \ge 1} \sum_{{\scriptstyle s^2+\sum\limits_{i=1}^{2s}id_i'=n,
  \atop\scriptstyle d_i'\ge 0 (1\le i\le 2s),}\atop\scriptstyle d_{2s}'\not = 0} 
  t^{2d_1'+\ldots + 2d_{2s}' +2s} (1-{1\over t^2})^{f(d_1')}\ldots
  t^{2d_{2s}'+1}(1-{1\over t^2})^{f(d_{2s}')}q^n     \\
&+\sum_{n \ge 1} \sum_{{\scriptstyle s^2+\sum\limits_{i=1}^{2s}id_i'=n, \atop
  \scriptstyle d_i'\ge 0 (1 \le i \le 2s),} \atop\scriptstyle d_{2s}' = 0} 
  t^{2d_1'+\ldots + 2d_{2s}' +2s }(1-{1\over t^2})^{f(d_1')}\ldots
  t^{2d_{2s-1}'+2}(1-{1\over t^2})^{f(d_{2s-1}')}(1+t)q^n  \\
=&1+\sum_{s\ge 1} \quad \sum_{d_i'\ge 0 (1\le i\le 2s), d_{2s}' \not = 0} 
    t^{2 \sum_{i=1}^{2s} i d_i'+ \sum_{i=1}^{2s} i} 
    (1-{1\over t^2})^{\sum_{i=1}^{2s}f(d_i')}q^{\sum_{i=1}^{2s}id_i'+s^2}     \\
 &+\sum_{s\ge 1}\quad \sum_{d_i'\ge 0 (1\le i\le 2s), d_{2s}' = 0} 
    t^{2 \sum_{i=1}^{2s} i d_i'+ \sum_{i=2}^{2s} i } 
    (1-{1\over t^2})^{\sum_{i=1}^{2s}f(d_i')}q^{\sum_{i=1}^{2s}id_i'+s^2}(1+t)     \\
=&1+\sum_{s\ge 1}\quad \sum_{d_i'\ge 0 (1 \le i\le 2s), d_{2s}'\not = 0}
    t^{2s^2+s}q^{s^2}(t^2q)^{\sum_{i=1}^{2s}id_i'}
    (1-{1\over t^2})^{\sum_{i=1}^{2s}f(d_i')}    \\
 &+ \biggl (\sum_{s\ge 1} \quad \sum_{d_i'\ge 0 (1\le i\le 2s), d_{2s}' = 0} 
    t^{2s^2+s}q^{s^2}(t^2q)^{\sum_{i=1}^{2s}id_i'}
    (1-{1\over t^2})^{\sum_{i=1}^{2s}f(d_i')}    \\
 &+\sum_{s\ge 1} \quad \sum_{d_i'\ge 0 (1\le i\le 2s-1)} 
    t^{2s^2+s-1}q^{s^2}(t^2q)^{\sum_{i=1}^{2s-1}id_i'}
    (1-{1\over t^2})^{\sum_{i=1}^{2s-1}f(d_i')}  \biggr )  \\
=&1+\sum_{s\ge 1}\quad \sum_{d_i'\ge 0 (1\le i\le 2s)} 
    t^{2s^2+s}q^{s^2}(t^2q)^{\sum_{i=1}^{2s}id_i'}
    (1-{1\over t^2})^{\sum_{i=1}^{2s}f(d_i')}    \\
 &+\sum_{s\ge 1}\quad \sum_{d_i'\ge 0 (1\le i\le 2s-1)} 
    t^{2s^2+s-1}q^{s^2}(t^2q)^{\sum_{i=1}^{2s-1}id_i'}
    (1-{1\over t^2})^{\sum_{i=1}^{2s-1}f(d_i')}.    \\
\endalign$$
Let $J$ be the set consisting of all the $j$ with $d_j' > 0$. Then 
$\sum_{n\ge 0}B_{0, n}(x, y)q^n $ equals 
$$\align
& 
1+\sum_{s\ge 1}t^{2s^2+s}q^{s^2}\sum_{J\subset \{1, \ldots, 2s\}}(1-{1\over t^2})^{|J|}
   \sum_{d_j'>0 (j\in J)} (t^2q)^{\sum_{j\in J} {jd_j'}}\\
 &+\sum_{s\ge 1}t^{2s^2+s-1}q^{s^2}\sum_{J\subset \{1, \ldots, 2s-1\}}(1-{1\over t^2})^{|J|}
   \sum_{d_j'>0 (j\in J)}(t^2q)^{\sum_{j\in J} {jd_j'}}\\
=&1+\sum_{s\ge 1}t^{2s^2+s}q^{s^2}\sum_{J\subset \{1, \ldots, 2s\}} (1-{1\over t^2})^{|J|}
    \prod_{j\in J}({1\over 1-(t^2q)^j}-1)\\
 &+\sum_{s\ge 1}t^{2s^2+s-1}q^{s^2}\sum_{J\subset \{1, \ldots, 2s-1\}} (1-{1\over t^2})^{|J|}
    \prod_{j\in J}({1\over 1-(t^2q)^j}-1)\\
=&1+\sum_{s\ge 1}t^{2s^2+s}q^{s^2}\prod^{2s}_{j=1}
    (1+(1-{1\over t^2}){(t^2q)^j\over 1-(t^2q)^j})
 +\sum_{s\ge 1}t^{2s^2+s-1}q^{s^2}\prod^{2s-1}_{j=1}
    (1+(1-{1\over t^2}){(t^2q)^j\over 1-(t^2q)^j})\\
=&\sum_{s\ge 0} (xy)^{2s^2+s}q^{s^2} 
    \prod^{2s}_{j=1} {1-(xy)^{2j-2}q^j \over 1-(xy)^{2j}q^j}
+\sum_{s\ge 1}(xy)^{2s^2+s-1} q^{s^2}
    \prod^{2s-1}_{j=1} {1-(xy)^{2j-2}q^j \over 1-(xy)^{2j}q^j}. \qed
\endalign$$

\par\noindent
{\it Remark 3.15}. In view of Yoshioka's results over finite fields
(the Remark 4.5 in \cite{Yos}), we think that the following is
a better closed formula for ${\W Z}_a(x, y, q)$:
$${\sum_{n \in \Zee} (xy)^{(2n+a)^2 - (2n+a) \over 2} q^{(2n+a)^2 \over 4}
\over q^{1 \over 12} (1 - xyq)} \cdot \prod_{d \ge 1} 
{1 - (xy)^{2d-1} q^d \over 1 - (xy)^{2d} q^d}.  \eqno (3.16)$$
For instance, we can verify that the lower degree terms
in (3.10) and (3.16) coincide by using MAPLE. However, 
we are unable to show that (3.10) and (3.16) are equal.


\Refs

\widestnumber\key{MMM}

\ref \key Bru \by R. Brussee \paper Stable bundles on blown up surfaces 
\jour Math. Z. \vol 205 \yr 1990 \pages 551--565 \endref

\ref \key Che \by J. Cheah \paper On the cohomology of Hilbert schemes 
of points \jour J. Alg.  Geom. \vol 5 \yr 1996 \pages 479-511
\endref

\ref \key D-K \by V.I. Danilov, A.G. Khovanskii \paper Newton polyhedra and an algorithm for computing Hodge-Deligne numbers \jour Math. USSR Izvestiya 
\vol 29 \pages 279-298 \yr 1987
\endref

\ref \key Del \by P. Deligne \paper Th\' eorie de Hodge III
\jour I.H.E.S. Publ. Math. \vol 44 \pages 5-77 \yr 1974
\endref

\ref \key Don \by S.K. Donaldson \paper Anti-self-dual Yang-Mills 
connections over complex algebraic surfaces and stable vector bundles
\jour Proc. Lond. Math. Soc. \vol 50\pages 1--26 \yr 1985\endref



\ref\key E-G
\by G. Ellingsrud, L. G\" ottsche 
\paper Variation of moduli spaces and Donaldson invariants under 
change of polarization
\jour J. reine angew. Math. \vol 467\pages 1-49\yr 1995
\endref

\ref \key F-M \by R. Friedman, J. W. Morgan \paper On the diffeomorphism 
types of certain algebraic surfaces II \jour J. Differ. Geom. \vol 27 
\pages 371-398 \yr 1988 \endref

\ref\key F-Q
\by  R. Friedman, Z. Qin
\paper Flips of moduli spaces and transition formulas for Donaldson 
polynomial invariants of rational surfaces
\jour Comm. Anal. Geom. \vol 3 \pages 11-83\yr 1995
\endref

\ref\key Ful \by W. Fulton
\paper Introduction to toric varieties
\jour Annals of Mathematics Studies
\vol 131
\publ Princeton University Press \publaddr Princeton\yr 1993
\endref

\ref \key Got  \by L. G\" ottsche 
\paper Change of polarization and Hodge numbers of moduli spaces 
of torsion free sheaves on surfaces
\jour Math. Z. \vol 223 \pages 247-260 \yr 1996 
\endref

\ref \key G-S \by L. G\" ottsche, W. Soergel
\paper Perverse sheaves and the cohomology of Hilbert schemes of 
smooth algebraic surfaces \jour Math. Ann. \vol 296 \pages 235-245 \yr 1993
\endref

\ref \key H-S \by H.J. Hoppe, H. Spindler
\paper Modulr{\" a}ume stabiler 2-B{\" u}ndel auf Regelfl{\" a}cher
\jour Math. Ann. \vol 249 \pages 127-140 \yr 1980
\endref

\ref \key LiJ \by J. Li \paper Algebraic geometric interpretation of 
Donaldson's polynomial invariants \jour J. Differ. Geom. \vol 37 
\pages 417--466 \yr 1993\endref

\ref \key {L-Q} \bysame
\paper On blowup formulae for the $S$-duality conjecture of Vafa and Witten
\jour Preprint
\endref

\ref \key {Qi1} \by Z. Qin \paper Stable rank-$2$ sheaves on blownup surfaces
\jour Unpublished \endref

\ref \key {Qi2}\bysame   
\paper  Moduli spaces of stable rank-$2$ bundles on ruled surfaces
\jour Invent. Math. \vol 110 \yr 1992 \pages 615-626
\endref

\ref \key Uhl \by K. Uhlenbeck \paper Removable singularity in 
Yang-Mills fields \jour Comm. Math. Phys. \vol 83 \pages 11-29 \yr 1982
\endref

\ref \key V-W \by C. Vafa, E. Witten \paper A strong coupling test of $S$-duality \jour Preprint
\endref

\ref \key Yos \by K. Yoshioka \paper The Betti numbers of the moduli space of stable sheaves of rank $2$ on $\Bbb P^2$ \jour J. reine angew. Math. \vol 453 \pages 193-220 \yr 1994
\endref




\endRefs
\enddocument